\documentclass[10pt,twocolumn,twoside]{IEEEtran}

\usepackage{amsfonts}
\usepackage{amsmath}
\usepackage{tikz}

\usepackage{graphicx}
\usepackage{subcaption}
\usepackage{sidecap}

\usepackage{mathrsfs}

\def\c{ \mathcal{C}}
\newcommand{\bea}{\begin{eqnarray}}
\newcommand{\eea}{\end{eqnarray}}

\ifCLASSINFOpdf
\else
\fi
\begin{document}
\title{Convexity of Solvability Set of\\ Power Distribution Networks}
%
%
%

\author{Anatoly~Dymarsky, and 
        Konstantin~Turitsyn,~\IEEEmembership{Member,~IEEE}\thanks{A. Dymarsky is with the Center for Energy Systems, Skolkovo Institute of Science and Technology, and University of Kentucky e-mail: a.dymarsky@skoltech.ru.}
\thanks{K. Turitsyn is with Massachusetts Institute of Technology e-mail:  turitsyn@mit.edu.}
}

\maketitle

\begin{abstract}
The solvability set of a power network -- the set of all power injection vectors for which the corresponding Power Flow equations admit a solution -- is central to power systems stability and security, as well as to the tightness of Optimal Power Flow relaxations. Whenever the solvability set is convex, this allows for substantial simplifications of various optimization and risk assessment algorithms. In this paper we focus on the solvability set of power distribution networks and prove convexity of the full solvability set (real and reactive powers) for tree homogeneous networks with the same $r/x$ ratio for all elements. We also show this result can not be improved:  once the network is not homogeneous, the convexity is immediately lost. It is nevertheless the case that if the network is almost homogeneous, a substantial practically-important part of the solvability set is still convex. Finally, we prove convexity of real solvability set (only real powers) for any tree network as well as for purely resistive networks with arbitrary topology. 
\end{abstract}


%
\IEEEpeerreviewmaketitle

\section{Introduction}
The solvability set of power networks, the set of all physically possible power injections, is the central object for power system voltage stability and security analysis. Whenever the operational point approaches the boundary of solvability, the power system may become unstable causing a blackout. In practice, physical constraints which follow directly from the Power Flow equations are often amended by various additional constraints, e.g.~on maximal nodal injection or limitations on line transmission. The set of all powers satisfying both physical and the additional constraints is usually denoted as feasibility set. In this paper we primarily focus on the solvability set $F$, which we will define as the set of all power injections for which the corresponding Power Flow equations admit a solution. 

The problem of the geometry of the power flow feasibility set has a long history with first discussions appearing as early as 1975 \cite{Galiana75, Schweppe75}. The non-convexity of the general solvability sets have been explicitly demonstrated and discussed for example in \cite{Hiskins, Bernie05, Makarov08}. 
Remarkably, optimization algorithms for OPF relying on convex relaxation often yield an exact result \cite{LL, Molzahn}. This observation lead to a host of activities focusing on understanding applicability of convex optimization methods toward OPF. 
In particular it was established that different formulations of convex optimization algorithms for OPF with constraints are exact for certain types of networks, e.g.~acyclic networks satisfying a set of realistic constraints \cite{LL, Molzahn, LL2, Radial1,Radial2, Low1,Low2,LowTutorial,QCQ}. In a related development certain geometrical properties of the feasibility sets for such networks were scrutinized  in \cite{Geometry1, Geometry2, Geometry3, Geometry4}. These latter works provide a geometric explanation for the exactness of convex relaxation by showing that although the feasibility set is not convex, the Pareto front in such cases coincides with the Pareto front of the convex hull of the feasibility set. Accordingly, when the utility function is non-decreasing, the optimal point would belong to the Pareto front and hence would be the same for the original problem and the convex relaxed one. 

In this paper we focus on establishing convexity of solvability set for several types of distribution networks. Whenever convexity is present this clearly simplifies associated OPF problems, allowing for an efficient solution for any convex utility function. This makes our work closely related to the previous studies on the subject. At the same time there is a substantial methodological difference between establishing convexity of the full feasibility set in our case and e.g.~establishing exactness of convex relaxation discussed in \cite{Geometry1,Geometry2,LowTutorial,QCQ}. Speaking geometrically, latter works considered non-convex feasibility set and focused on a particular area of the boundary, the so-called Pareto front. This approach provides additional flexibility to incorporate various constraints, but it is not informative about the internal points of $F$ and can only explore the geometry of $F$ locally. Our approach surveys the whole boundary of the solvability set to detect potential boundary non-convexity. Crucially, it also involves a topological argument of \cite{Dymarsky1, Dymarsky2} which connects geometry of the boundary and convexity of the interior by ensuring that $F$ has no ``holes" inside.  Thus our approach is complementary to the convex optimization-based techniques and can establish global geometric properties of $F$, beyond the scope of the traditional convex analysis. 

This paper is organized as follows. In the next section we introduce the notations and formulate the sufficient condition for convexity of \cite{Dymarsky2}. We also explain its geometric interpretation. In section \ref{sec3} we prove convexity of the full solvability set for homogeneous distribution networks with the tree topology. In section \ref{realsetunbalanced} we discuss convexity of the real solvability set and prove its convexity for several different types of networks: arbitrary distribution networks with the tree topology and purely resistive networks with arbitrary topology. We conclude in section \ref{conclusions}.

\section{Notations and methodology}
\subsection{Model and Notations}
Throughout this paper we consider an AC model of power system consisting of  $n$ PQ-buses $i=1,\dots,n$ and a slack bus $i=0$. The network is parametrized by a symmetric complex-valued admittance matrix $Y_{ik}$, $i,k=0,\dots,n$,
\begin{eqnarray}
Y_{ik}&=&\left\{
\begin{array}{c r}
\sum_{l\sim i} y_{il} & {\rm \ if\ }i=k\\
-y_{ik} & {\rm \ if\ }i\sim k\\
0 & {\rm \ if\ }i\not\sim k\\
\end{array}
\right.
\end{eqnarray}
By $i \sim k$ we denote that nodes $i$ and $k$ are connected. We further assume that in that case 
$y_{ik}\neq 0$.
The power injection at node $i$ is given by (for $i=1,\dots,n$),
\begin{equation}
\label{PF}
S_i=P_i+jQ_i=V_i \sum_{k\sim i} y^*_{ik}(V_i-V_k)^* \ .
\end{equation}
In our notations negative $P_i$ corresponds to power consumption and positive $P_i$ to power injection at node $i$. It is convenient to combine real and reactive powers into a vector $p=(P_1,\dots,Q_1,\dots)^T\in {\mathbb R}^{2n}$. The complex voltages $V_i$ for $i=1,\dots,n$ can be combined into  a complex vector ${\bf V}=(V_1,\dots)^T\in {\mathbb C}^n$ while slack voltage is taken to be $V_0\equiv 1$. The full solvability set $F\subset {\mathbb R}^{2n}$ is a set of all points $p\in {\mathbb R}^{2n}$ such that the system of Power Flow equations (\ref{PF}) is feasible (have a solution). Similarly we introduce the solvability set of real powers, or real solvability set, as the combination of all points $(P_1,\dots,P_n)^T\in {\mathbb R}^n$ such that the system of equations (\ref{PF}) is feasible at least for some $Q_i$.


One of the central objects in our approach is a linear combination of real and reactive powers  
\begin{eqnarray}
\label{cp}
c\cdot p=\sum_{i=1}^n c_i P_i+ \sum_{i=1}^n {c}_{n+i} Q_i &=&\sum_{i=1}^n\Re({\mathcal C}^*_i S_i)\ ,\\
\mathcal{C}_i&=&{c}_i+j{c}_{n+i}\ .
\label{cdef}
\end{eqnarray}
Here we introduced $n$ complex variables ${\mathcal C}_i$ which encode a vector $c\in {\mathbb R}^{2n}$. For each $c\neq 0$, minimization problem 
\begin{equation}
\label{supportingH}
\arg \min_{p\in F} c\cdot p\ ,
\end{equation}
gives the intersection points of the solvability set $F\subset {\mathbb R}^{2n}$ 
and the supporting hyperplane orthogonal to $c$. The same linear combination can be rewritten in matrix notations as
\begin{eqnarray}
\label{cp2}
c\cdot p&=&{\bf V}^\dagger H^{}_c {\bf V}-{\bf V}^\dagger J^{}_c -J_c^\dagger {\bf V}\ ,\\
\label{Hc}
(H^{}_c)^{}_{ik}&=&\left\{
\begin{array}{l r}
\Re({\mathcal C}_i \sum_{l\sim i} y_{il}) & {\rm \ if\ }i=k\\
-(\c_i y_{ik}+\c_k^* y_{ik}^*)/2 & {\rm \ if\ }i\sim k\\
0 & {\rm \ if\ }i\not\sim k\\
\end{array}
\right.\\
\label{Jc}
(J_c)_i&=&\left\{\begin{array}{l r}
 \c_i y_{i0}/2 & \qquad \qquad \ \  \quad {\rm if\ }i\sim 0 \\
0  &  {\rm if\ } i\not\sim 0
\end{array}\right.\
\end{eqnarray}

\subsection{Sufficient Condition for  Convexity}
Real and reactive powers are the quadratic functions of voltages, $p=p(V)$. Hence the solvability set is an image of a particular quadratic map defined by (\ref{PF}). 
To establish convexity of solvability set we will rely on the result (Proposition 2') of \cite{Dymarsky2}, which proves the sufficient condition for an image of a quadratic map to be convex. Below we formulate that sufficient condition and explain the geometric intuition behind it. 

First of all, the sufficient condition requires the corresponding quadratic map to be {\it definite}, i.e.~there must exist a vector $c_+\neq 0$ such that the quadratic function $c_+\cdot p(V)$ is bounded from below.  Geometrically, it means that the corresponding solvability set is confined to a half-space defined by an appropriate hyperplane.  With help of (\ref{cp2}) this can be reformulated as positive definiteness of $H_+\equiv H_c$ (\ref{Hc}) for some  $c=c_+$. If all lines have non-vanishing resistance, $\Re(y_{ik})>0$, the corresponding vector $c_+$ is readily given by 
\bea
\label{cplus}
(c_+)_i=\left\{
\begin{array}{l r}
1, & 1\le i\le n \\
0, & i>n
\end{array}
\right.
\eea
This is simply the observation that in a network with resistive links the total power consumption $-(P_1+\dots+P_n)$ is limited from above. To see that we write 
\bea
{\bf V}^\dagger H_+ {\bf V}=\sum_{i\sim k, i>k>0} \Re(y_{ik})|V_i-V_k|^2+\sum_{i\sim0} \Re(y_{i0})|V_i|^2, 
\label{Hplus}
\eea
which is manifestly positive unless all $V_i=0$. We emphasize that (\ref{Hplus}) establishes {\it definiteness} of quadratic maps associated with both full solvability set and the solvability set of only real powers.  

The main idea of the geometric approach of \cite{Dymarsky2} is to focus on definite maps and to study boundary points  lying at an intersection of the solvability set with a supporting hyperplane defined by some vector $c\neq 0$. These points are given by (\ref{supportingH}). Minimization problem (\ref{supportingH}) has a solution provided $H_c\succeq 0$ and the following equation is feasible 
\bea
\label{boundary}
H_c {\bf V}_{\rm b}=J_c\ .
\eea
When $H_c\succ 0$ the solution alway exists and unique, but when $H_c$ is singular we have Fredholm alternative: the solution either does not exist or it is not unique because there are non-zero vectors  ${\bf V}_{\rm null}$ such that 
\bea
\label{flatedge}
H_c {\bf V}_{\rm null}&=&0\ ,\\
J_c^\dagger {\bf V}_{\rm null}&=&0\ .
\label{flatedge2}
\eea
In the latter case linear combination of voltages ${\bf V}={\bf V}_{\rm b}+{\bf V}_{\rm null}$ with any ${\bf V}_{\rm null}$ satisfying (\ref{flatedge},\ref{flatedge2}) gives boundary points $p(V)$ which minimize (\ref{supportingH}) and hence belong to a ``flat edge" on the boundary of $F$. This boundary region is typically non-convex \cite{Dymarsky2}. On the contrary, if for any $H_c\succeq 0$ the system of equations (\ref{boundary},\ref{flatedge},\ref{flatedge2}) admits no solutions, the full solvability set is convex \cite{Dymarsky2}. 

\noindent{\bf Sufficient condition for convexity (I).} To summarize,  Proposition 2' of \cite{Dymarsky2} ensures convexity of solvability set provided that: (i) there exist vector $c_+$ such that $H_+\succ 0$ and (ii) for any $c\neq 0$ such that $H_c\succeq 0$ and singular equation (\ref{boundary}) has no solutions. 

The same condition can be conveniently reformulated in the following way. We introduce a $(n+1)\times (n+1)$ Hermitian matrix
\bea
\label{matrixA}
A(a,c)=\left(
\begin{array}{c|c}
a & -J_c^\dagger \\
\hline\\[-8pt]
-J^{}_c & H^{}_{c}
\end{array}
\right)
\eea
and notice that (\ref{boundary}) together with $H_c\succeq 0$ is equivalent to  matrix $A$ with $a={\bf V}_{\rm b}^\dagger H_c {\bf V}_{\rm b}\ge 0$ being positive semi-definite and singular with the eigenvector associated with zero eigenvalue given by $(V_0,{\bf V}^T_{\rm b})^T$. Furthermore the feasibility of (\ref{flatedge},\ref{flatedge2}) means that zero eigenvalue of $A$ is at least double degenerate. \\
\noindent{\bf Sufficient condition for convexity (II).} The convexity condition of \cite{Dymarsky2} can be reformulated as follows. If the corresponding quadratic map is definite and for any $c\neq 0$ and $a$ such that $A\succeq 0$, matrix $A$ can have at most one zero eigenvalue, the associated solvability set is convex.

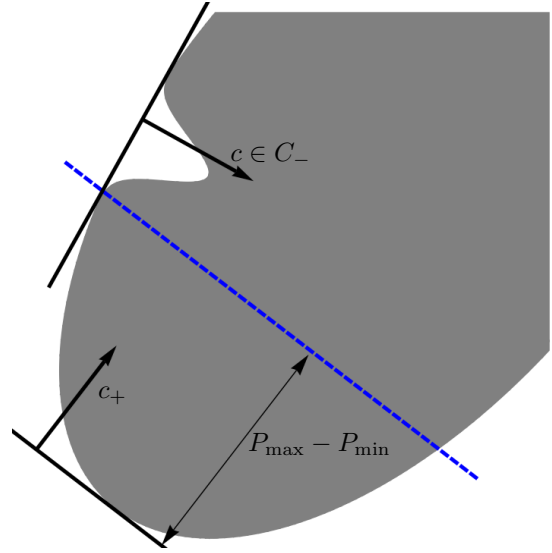
\begin{figure}[t]
\begin{tikzpicture}
\node[inner sep=0pt] (picture) at (3.5,3.4)
    {\includegraphics[width=0.4\textwidth]{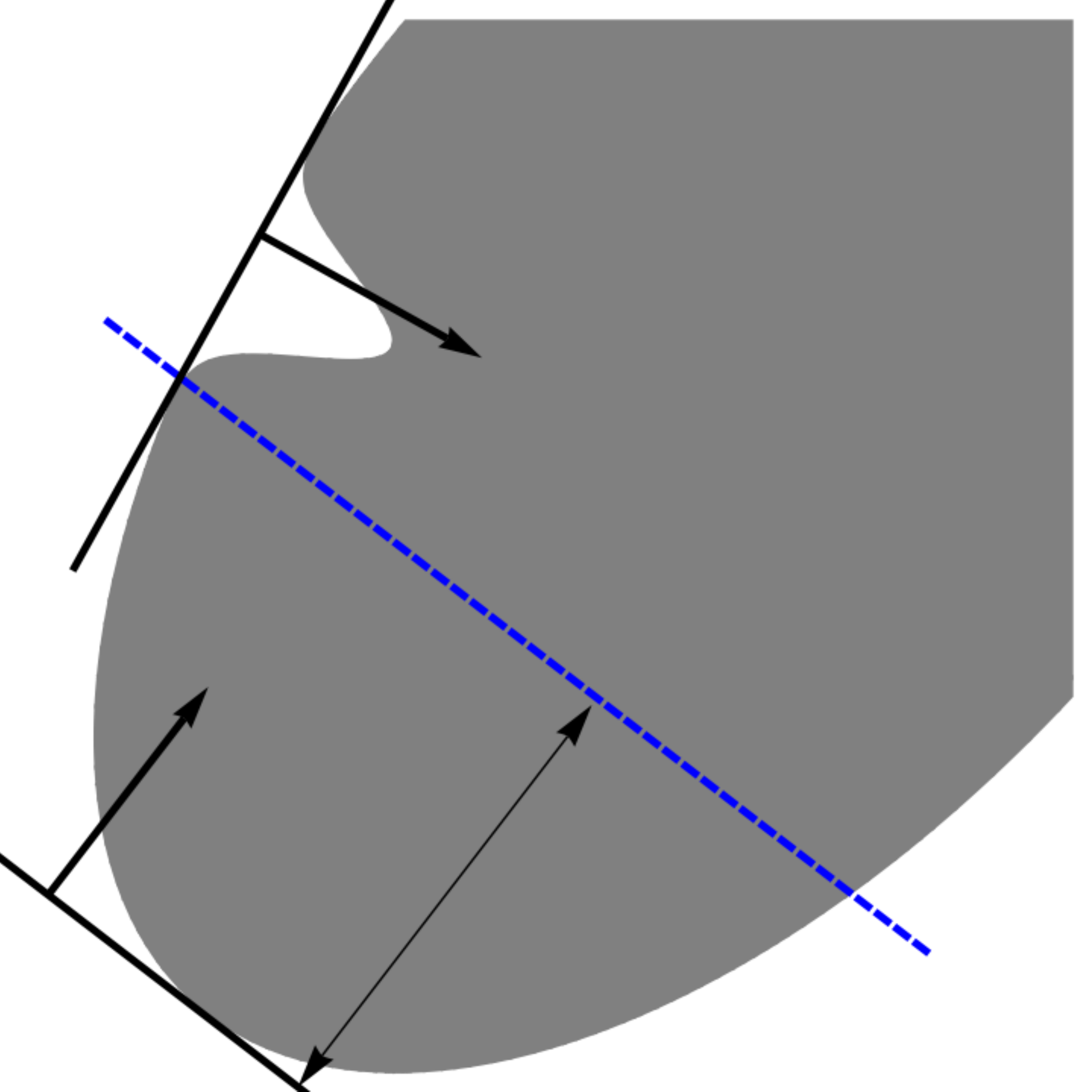}};
                \node at (1.2,1.8) {$c_+$};
                \node at (3.95,1.15) {$P_{\rm max}-P_{\rm min}$};    
                \node at (3.3,5) {$c\in C_-$};  
\end{tikzpicture}
\caption{The idea behind identifying convex subset within the solvability set: to identify maximal $P_{\rm max}$ such that 
all boundary non-convenvexities satisfy $c_+\cdot p\ge P_{\rm max}$. Then the subset of the solvability set defined by $P_{\rm min}\le c_+\cdot p\le P_{\rm max}$ is convex.}
\label{fig_sim}
\end{figure}

\noindent{\bf Convex subset.}
Whenever for some vector $c\neq 0$ matrix $H_c\succeq 0$ and the set of equations (\ref{boundary},\ref{flatedge},\ref{flatedge2}) admit a non-trivial solution with ${\bf V}_{\rm null}\neq 0$, the full solvability set is likely not to be convex. Let us denote a set of all such vectors $c\neq 0$  as $C_-$. For a given $c\in C_-$ the potentially non-convex  boundary region can not stretch beyond a particular hyperplane defined by an equation $c_+ \cdot p=P(c)$, where 
\bea
\label{P(c)}
P(c)&=&\min_{{\bf V}_{\rm null}} {\bf V}^\dagger\, H_+ {\bf V}-J^\dagger_+ {\bf V}-{\bf V}^\dagger J_+\ ,\\
&&{\rm where\ \ \ }{\bf V}={\bf V}_{\rm b}+{\bf V}_{\rm null}\ . \nonumber
\eea 
Here minimization goes over all vectors ${\bf V}_{\rm null}$ satisfying  (\ref{flatedge},\ref{flatedge2}). Now, if by $P_{\rm max}$ we denote minimum of $P(c)$ for all possible vectors $c\in C_-$, then the subset of $F$ defined by an inequality $c_+\cdot p\le P_{\rm max}$ will be convex (Proposition 2' of \cite{Dymarsky2}). The geometrical intuition here is that as soon as one can make sure all non-convex boundary regions do not stretch beyond i.e.~do not intersect a particular hyperplane, the compact subset of $F$ constrained by that hyperplane is convex. We can also introduce $P_{\rm min}$ as the global minimum of $c_+\cdot p$, such that for any point $p$ from the solvability set $c_+\cdot p \ge P_{\rm min}$. The the convex subset can be conveniently define by the inequalities $P_{\rm min}\le c_+\cdot p\le P_{\rm max}$.
This is illustrated in Fig.~\ref{fig_sim}.

\section{Convexity of Full Solvability Set for Homogeneous Distribution Networks}
\label{sec3}
In this section we prove convexity of the full solvability set for homogeneous distribution networks. We would say the network is homogeneous if all $y_{ik}$ for any pair of connected nodes $i\sim k$ have the same argument $\arg y_{ik}=\phi$. General linear transformation does not affect convexity. It is therefore convenient to perform the following linear change of variables $P_i,Q_i$, and $c_i$
\bea
\label{phi}
S_i\rightarrow e^{i\phi}S_i,\quad \c_i\rightarrow e^{i\phi}\c_i\ ,
\eea
amended by the redefinition $y_{ik}\rightarrow e^{-i\phi} y_{ik}$. This change leaves
the Power Flow equations (\ref{PF}) as well as the equations (\ref{cp},\ref{cp2}) invariant.  In the new variables the network is purely resistive with $y_{ik}$  for any two connected nodes $i\sim k$ being positive real numbers. We furthermore assume that the network is connected, i.e.~any two nodes can be connected through a combination of links. 

\noindent{\bf Lemma 1.} For a purely resistive connected network and vector $c\neq 0$ such that $H_c\succeq 0$, either $\Re(\c_i)>0$ for all $i=1,\dots,n$ or  $\c_i=j \alpha$  where $\alpha$ is the same real number for all $i=1,\dots,n$ and $H_c=0$.\\
{\bf Proof.} The condition $H_c\succeq 0$ requires all diagonal elements of $H_c$ to be non-negative. Since all $y_{ik}$ are positive, from (\ref{Hc}) it follows that  $\Re(\c_i)\ge0$ for all $i=1,\dots,n$. Now let us assume that  $\Re(\c_i)=0$ for a particular $i$. For $H_c$ to be positive semi-definite all elements $(H^{}_c)^{}_{ik}$  must be zero, $(\c_i+\c_k^*)y_{ik}=0$. From here it follows that $\Re(\c_k)=0$ and $\Im(\c_k)=\Im(\c_i)\equiv\alpha$ for all $k\sim i$. And since the network is connected, repeating this consideration enough times, we find $\Re(\c_k)=0$ and $\Im(\c_k)=\alpha$ for all $k=1,\dots,n$. 

\subsection{Convexity of homogeneous acyclic networks}
\label{balancedtree}
First we prove convexity of the full solvability set for homogeneous networks with tree topology. The idea of the proof is to use Theorem 3.4 of \cite{Holst} to show that a positive semi-definite matrix $A$ (\ref{matrixA}) can have at most one zero eigenvalue. To that end it is enough to show that $A$, when positive semi-definite, has the topology of the network graph. This is almost apparent from the form of $H_c$ and $J_c$ but one needs to prove that $(H_c^{})^{}_{ik}$ and $(J_c^{})^{}_i$ do not vanish when $i \sim k$ and $i \sim 0$ correspondingly. 

\noindent{\bf Theorem 1.} The full solvability set of a homogeneous connected acyclic (with tree topology) network is convex.\\ 
{\bf Proof.} By applying the linear transformation \eqref{phi} a homogeneous network can be brought to the form of a purely resistive one. As follows from the explicit form of \eqref{Hplus} matrix $H_c$ with $c=c_+$ given by \eqref{cplus} is positive-definite, which establishes that the associated quadratic map is definite. Next, we consider vector $c\neq 0$ and real number $a$ such that matrix $A$ given by \eqref{matrixA} is positive semi-definite. This in particular implies that any sub-matrix of $A$ is also positive semi-definite. Hence $H_c\succeq 0$ and from Lemma 1 it follows that either $H_c=0$ or $\Re(\c_i)>0$ for all $i=1,\dots,n$. First we consider the former possibility. For $A$ to be positive semi-definite while $H_c=0$,  it requires  $J_c=0$. Hence for all $i\sim0$, $\c_i y_{i0}=0$. Since  all $\c_i=j\alpha$ this immediately implies all $\c_i=0$, which contradicts the assumption that $c\neq 0$. Hence, from $A\succeq 0$ and $c\neq 0$  it follows that $\Re(\c_i)>0$ for all $i=1,\dots,n$. Consequently $\Re((H^{}_c)^{}_{ik})=-\Re(\c_i+\c_k)y_{ik}>0$ is non-zero for $i\sim k$, and  $\Re((J_c)_{i}^{})=-\Re(c_i)y_{i0}>0$ is non-zero for $i\sim 0$. This is sufficient to use Theorem 3.4 of \cite{Holst} to establish that $A$ has at most one zero eigenvalue.  Now the sufficient condition for convexity (II) applies, which completes the proof.

\subsection{Radial 3-bus model}
\label{3bus}
As a next step we demonstrate that Theorem 1 can not be strengthen. Namely, introducing non-zero relative phases in the admittance coefficients may immediately render the solvability set non-convex. To show that we consider a simple radial 3-bus model with two PQ-buses connected in series to a slack bus, resulting in the following admittance matrix 
\begin{eqnarray}
Y=\left(
\begin{array}{ccc}
y & -y & 0\\
-y & (y+1) & -1 \\
0& -1 & 1
\end{array}
\right)\ .
\end{eqnarray}
Here without loss of generality we chose one admittance to be a complex parameter $y_{01}=y$, while the other is $y_{12}=1$. We will see below that unless $\Im(y)=0$ there is a non-empty set $C_-$ and vector(s) $c\in C_-$ such that $H_c\succeq 0$  and \eqref{boundary} is feasible together with (\ref{flatedge},\ref{flatedge2}) and a non-zero ${\bf V}_{\rm null}$. Let us assume $c\in C_-$ and also $\c_1=c_1+j c_3$ is non-zero. Using $J_c=(\c_1 y/2,0)^T$ and \eqref{flatedge2} we find vector  ${\bf V}_{\rm null}$ up to an overall multiplier, ${\bf V}^\dagger_{\rm null}=(0,1)$. For this vector to be annihilated by the matrix 
\begin{eqnarray}
H_c=\left(
\begin{array}{cc}
\Re(\c_1 y+\c_1) & -(\c_1 +\c_2^*)/2 \\
-(\c^*_1+\c_2)/2 & \Re(\c_2)
\end{array}
\right)
\end{eqnarray}
we must require $\Re(\c_2)=0$ and $(\c_1 +\c_2^*)=0$. 
This readily gives that $\c=i(1,1)^T$ up to an overall real coefficient and ${\bf V}_{\rm b}=(-i y/(2\Im y),0)^T$. That $C_-$ is non-empty indicates solvability set may have a  non-convex boundary spanned by the points 
\bea
\label{noncb}
S_i(z)=S_i({\bf V}_{\rm b}+z{\bf V}_{\rm null})
\eea for arbitrary complex $z$. A straightforward check confirms that.  Indeed all points \eqref{noncb} lie on a hyperplane 
\bea
\label{edgehp}
\Im(S_1(z)+S_2(z))=-{|y|^2 \over 4\Im(y)}\ ,
\eea
and these are the only points of the solvability set that belong to this hyperplane. Equation \eqref{noncb} defines a parabolic surface inside \eqref{edgehp}, which is obviously non-convex. That establishes non-convexity of the full solvability set. 

It is important to note that $C_-$ is non-empty and boundary non-convexity exist whenever $\Im y\neq 0$ no matter how small it is. This result shows that even in a simplest network with tree topology full solvability set ceases to be convex unless all admittance coefficients $y_{ik}$ have the same complex phase. 

This result might seem surprising at first. It can be shown that the sufficient condition of \cite{Dymarsky2} establishes strong convexity of the full  solvability set when the network is homogeneous. Naively strong convexity is stable i.e.~it can not be destroyed by an infinitesimal change of parameters. The caveat here is that the full solvability set is non-compact. Hence the non-convexity appearing whenever the network is non-homogeneous is located very far from the origin. This can be readily seen in the 3-bus example: when $\Im (y)\neq 0$ but very small the non-convexity is located very far away from the operationally important region. Since $C_-$ only includes vectors collinear to $\c=i(1,1)^T$ it is easy to calculate \eqref{P(c)}
\bea
P(c)&=&\min_z \Re(S_1(z)+S_2(z))=\\
&&\min_z \Re(y)\left|z+{iy\over 2\Im y}\right|^2+{(\Re y) |y|^2\over 4(\Im y)^2}\ ,\nonumber
\eea
and 
\bea
\label{Pmax}
P_{\rm max}={(\Re y) |y|^2\over 4(\Im y)^2}\ .
\eea
Thus we find that the subset of the full solvability set constrained by a linear inequality $P_1+P_2\leq P_{\rm max}$ is convex.  When $\Im y$ is small \eqref{Pmax} guarantees that a very large subset of the full solvability set is convex.

Finally, we note that the total power consumption $-(P_1+P_2)$ is bounded from above by $-P_{\rm min}={|y|^2\over 4\Re(y)}$, hence the feasible regimes confined to the convex subset of the full solvability set are 
\bea
\label{Pmaxmin}
-{|y|^2\over 4\Re(y)}=P_{\rm min}\le P_1+P_2\le P_{\rm max}={(\Re y) |y|^2\over 4(\Im y)^2}\ .
\eea

\begin{figure*}[!t]
\centering
\begin{subfigure}[t]{0.3\textwidth}
        \centering
        \includegraphics[width=\textwidth]{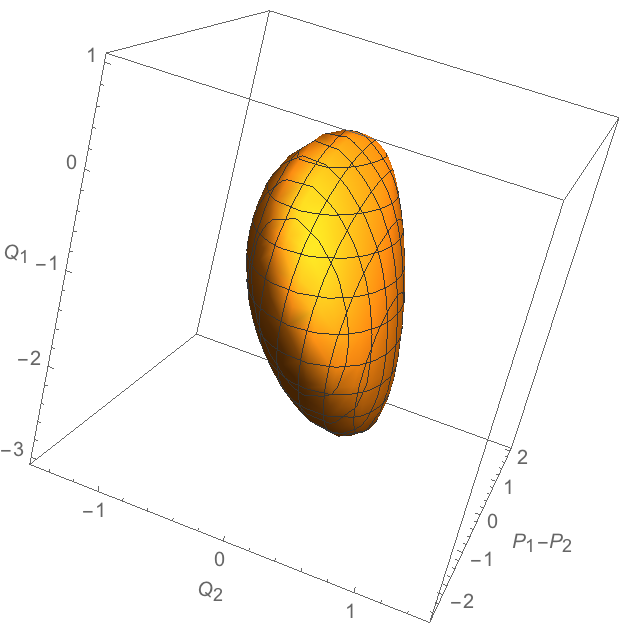}
        \caption{}
        \label{fig1}
    \end{subfigure}
  \hfill
    \begin{subfigure}[t]{0.3\textwidth}
        \centering
        \includegraphics[width=\textwidth]{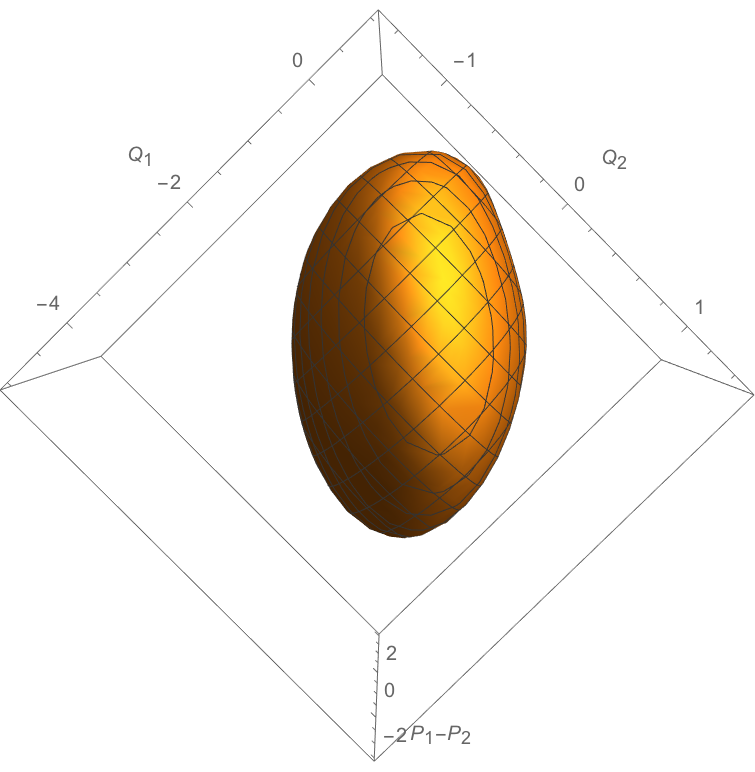}
        \caption{}
        \label{fig2}
    \end{subfigure}
 \hfill
    \begin{subfigure}[t]{0.3\textwidth}
        \centering
        \includegraphics[width=\textwidth]{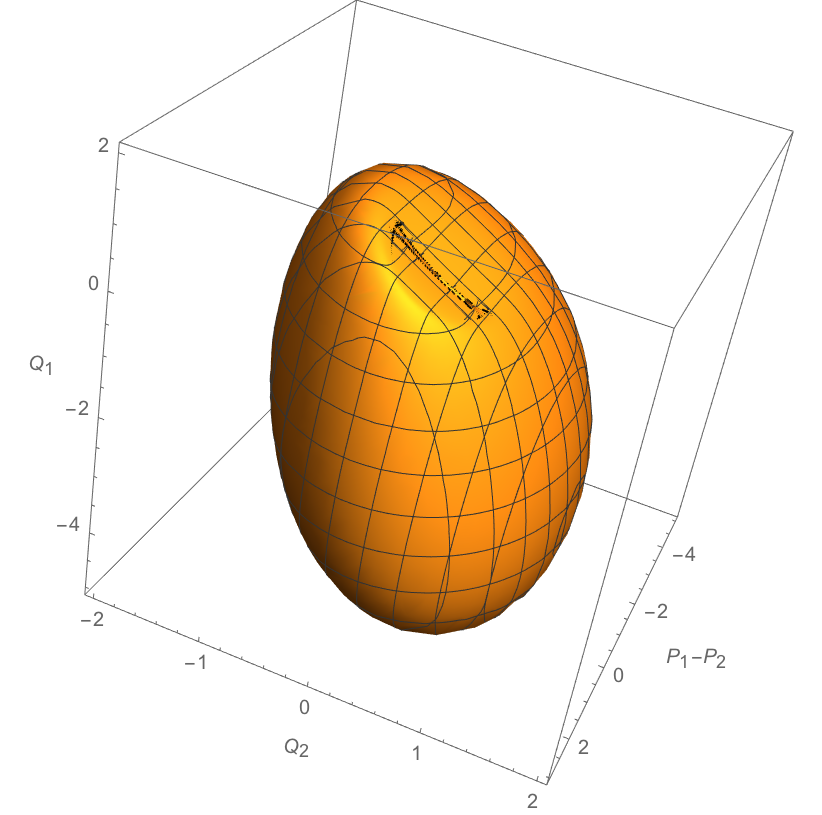}
        \caption{}
        \label{fig3}
    \end{subfigure}
\caption{Projections of the full solvability set for the 3-bus model discussed in section \ref{3bus} on the hyperplanes $P_1+P_2={\rm const}$. (a)  $P_1+P_2=0$ and the projection is strongly convex. (b) $P_1+P_2=1/2$ and the projection develops a flat boundary i.e.~it is convex but not strongly convex. (c) $P_1+P_2=1$ and
 the projection is visibly non-convex.}
\label{fig_triptih}
\end{figure*}

We show the projections of the full solvability set on different hyperplanes $P_1+P_2={\rm const}$ for $y=1+j$ in Fig.~\ref{fig_triptih}. In this case $P_{\rm min}=-1/2$ and $P_{\rm max}=1/2$. For $P_1+P_2=0<P_{\rm max}$ the projection is strongly convex, Fig.~\ref{fig_triptih} (a). The value $P_1+P_2=1/2=P_{\rm max}$ is critical as in this case the projection of the solvability set develops a flat edge, Fig.~\ref{fig_triptih} (b). Finally, for $P_1+P_2=1>P_{\rm max}$, it is non-convex, Fig.~\ref{fig_triptih} (c).

\section{Convexity of Real Solvability Set}
\label{realsetunbalanced}
As we saw above the full solvability set of the unhomogeneous networks is generally non-convex. Nevertheless convexity often can be preserved if instead of the full solvability set one considers only the real solvability set -- the solvability set of real powers. From the point of view of the underlying quadratic map, we restrict our consideration to $n$ real-valued quadratic functions, $P_i=\Re(S_i)$, \eqref{PF}  of $n$ complex variables $V_i$. Hence  in  \eqref{supportingH} we should take $m=n$ with $p,c\in {\mathbb R}^n$. In practice it is convenient to work with the full vectors $p=(P_1,\dots,Q_1,\dots)^T\in {\mathbb R}^{2n}$ and $c\in {\mathbb R}^{2n}$ but restrict $c_i=0$ for $i>n$. The same condition can be written as $\Im(\c_i)=0$. Notice, that \eqref{cplus} satisfies this condition and therefore if all lines have non-zero resistance the associated quadratic map is definite. In what follows we will assume the latter, i.e.~$\Re(y_{ik})>0$ for all $i\sim k$. In case  some lines are purely reactive with $\Re(y_{ik})=0$, it is possible to add a small positive number $\epsilon$ to each such $y_{ik}$, establish convexity in that case and then take $\epsilon$ to zero. Since the convexity is stable property, this would grantee convexity for $\epsilon=0$ as well.

\subsection{Convexity of acyclic networks}
If the consideration is restricted to real powers, real solvability set of any radial (acyclic) network is convex. As in the case of full solvability set for homogeneous acyclic networks, section \ref{balancedtree}, the idea of the proof would be to use Theorem 3.4 of \cite{Holst}.\\
\noindent{\bf Lemma 2.} For a connected network with non-zero resistances $\Re(y_{ik})>0$ for all $i\sim k$, and vector $c\neq 0$, $\Im(\c_i)=0$, such that $H_c\succeq 0$ all $\c_i>0$.\\
{\bf Proof.} The proof is very similar to the proof of Lemma 1. The condition $H_c\succeq 0$ requires all diagonal elements of $H_c$ to be non-negative. Since all $\Re(y_{ik})$ are positive and $\c_i$ are real and can be brought outside of $\Re$, from (\ref{Hc}) it follows that  $\c_i\ge0$ for all $i=1,\dots,n$. Now let us assume that  $\c_i=0$ for a particular $i$. For $H_c$ to be positive semi-definite all elements $(H^{}_c)^{}_{ik}$  must be zero, $(\c_i+\c_k)y_{ik}=0$. This is only possible if $\c_k=0$ for all $k\sim i$. And since the network is connected, repeating this consideration enough times, we find $\c_k=0$ for all $k=1,\dots,n$. This contradicts the assumption $c \neq 0$, and therefore all $\c_i>0$.

\noindent{\bf Theorem 2.} The real solvability set of any connected acyclic (with tree topology) network is convex.\\ 
{\bf Proof.} As was explained in the beginning of this section, without loss of generality we can assume that $\Re(y_{ik})>0$ for any two connected nodes $i\sim k$. Then \eqref{cplus} and \eqref{Hplus} establish that the corresponding quadratic map is convex. To apply the sufficient condition for convexity (I) we assume that number $a$
and $c\neq 0$, $\Im(\c_i)=0$ are such that matrix $A$ of \eqref{matrixA} is positive semi-definite. Using Lemma 2 for the submatrix $H_c$ of $A$ we readily find that all $\c_i$ are positive. From here it follows that matrix elements $(H_c)_{ik}$ for any two connected nodes $i\sim k$ are non-zero, 
\bea
\Re(H^{}_c)_{ik}=-(\c_i+\c_j)\Re(y^{}_{ik})<0\ .
\eea
Similarly for $i\sim 0$ vector element $(J_c)_i$ also does not vanish,
\bea
\Re((J_c)_i)=\c_i \Re(y_{i0})>0\ .
\eea
Hence, matrix $A$ (\ref{matrixA}) has the topology of the network. For acyclic networks we can immediately apply Theorem 3.4 of \cite{Holst} which finishes the proof. 

\subsection{Convexity of purely resistive networks}
When the network is purely resistive $\Im(y_{ik})=0$ and $\Re(y_{ik})>0$ for $i\sim k$ the real solvability set is also convex for any network topology. The approach we take here to prove it is similar to the proof in the subsection 5.1 of \cite{Dymarsky2} that the power solvability region for purely resistive DC model is convex. Notice however that in this paper we do not restrict voltages to be real, i.e.~we consider an AC model of power flow with purely resistive connected network of arbitrary topology.\\ 
\noindent{\bf Theorem 3.} The real solvability set of any connected purely resistive  network is convex.\\ 
{\bf Proof.} The idea of the proof is to use the sufficient condition for convexity (I). We start by using \eqref{cplus} and \eqref{Hplus} to show that the associated quadratic map is definite. Next, we consider $c\neq 0$ with all $\c_i$ being real $\c_i$ such that $H_c\succeq 0$ and singular. Using Lemma 2 we find that all $\c_i$ are strictly positive. Let us consider a normalized vector $|{\bf V}_{\rm null}|=1$, ${\bf V}_{\rm null}\in {\mathbb C}^n$ which satisfies \eqref{flatedge}. This vector minimizes the  following quadratic form 
\bea
\label{QF}
{\bf V}^\dagger H_c {\bf V}= \sum_{i=1}^n(H_c)_{ii} |V_i|^2 -2\sum_{i\sim k} (\c_i+\c_j)y_{ik} \Re(V_i^*V_k)\ ,
\eea
subject to ${\bf V}^\dagger {\bf V}=1$. Taking into account that for any $i\sim k$ the combination $ (\c_i+\c_j)y_{ik}>0$, for the given values of $|V_i|$, 
the quadratic form (\ref{QF}) will be minimal if $\Re(V_i^*V_k)=|V_i||V_k|$. Hence all components of ${\bf V}_{\rm null}$ minimizing (\ref{QF}) must have  the same phase. Thus, without loss of generality we can assume that  all components of  ${\bf V}_{\rm null}$ are real and non-negative. 
Therefore 
\bea
J_c^\dagger {\bf V}_{\rm null}=\sum_{i\sim 0} \c_i y_{i0} V_i\ , 
\eea
is a sum of non-negative terms. For the condition (\ref{flatedge2}) to be satisfied, it would require $V_i=0$ for each node connected with the slack $i\sim 0$.  Next, we  write $i$-th component of the vector $H_c {\bf V}_{\rm null}=0$, 
\bea
\c_i\sum_{k\sim i}y_{ik} V_i-\sum_{k\sim i}(\c_i+\c_k)y_{ik} V_k/2=0\ .
\eea
Assuming $V_i=0$, since all $V_k$ are non-negative, the condition \eqref{flatedge} can be only satisfied if $V_k=0$ at all nosed $k$ adjacent to $i$. Since the network is connected, by repeating this logic we find that all components of vector ${\bf V}_{\rm null}$ are zero, which contradicts the assumption that ${\bf V}_{\rm null}$ is normalized. Finally, we conclude that whenever $H_c\succeq 0$ and $c\neq 0$, the system of equations(\ref{flatedge},\ref{flatedge2}) is not feasible, the sufficient condition for convexity (I) applies and the real solvability set is convex.

\section{Conclusion}
\label{conclusions}
In this paper we have employed a novel geometric approach of \cite{Dymarsky2} to establish convexity of power flow solvability sets for various typos of networks. For the distribution network with a radial topology, we established that the full solvability set consisting of real and reactive powers is convex whenever the network is homogeneous (Theorem 1). By considering a simple example of a radial 3-bus network we have demonstrated that the condition of the network to be homogeneous is crucial, once it is relaxed the full solvability set immediately ceases to be convex. Next, we considered the real solvability set, i.e.~the solvability set consisting only of real powers, and established its convexity for an arbitrary distribution network with tree topology (Theorem 2). We have also shown that real solvability set is convex for a purely resistive network with arbitrary topology (Theorem 3).

Our results open a few new directions for research. First, it remains an open question to identify other network classes with convex solvability sets. For example, our preliminary studies suggest that a homogeneous distribution network with a ring topology ($n$ PQ-buses connected in a series with the first and the last connected to the slack) might have a convex full solvability set. Second, the approach of \cite{Dymarsky2} allows to establish rigorous analytic bounds on the convex subset within the solvability set. Realistic networks may not satisfy all necessary conditions for the solvability set to be convex. Thus, transmission networks are often close to be homogeneous, but variations of the impedance to resistance ratio are always present. As was illustrated in section \ref{3bus} this is still likely to be sufficient for the practically-relevant part of the solvability set to be convex. We leave it as an interesting problem for the future to derive analytic bounds similar to (\ref{Pmax}) for various types of networks. 

Certification of solvability set convexity for realistic power grid would allow for application of many powerful algorithms available for convex domains. For example, for the purposes of  robust power flow analysis, construction of inner approximations of secure operation regions could be accomplished by taking a convex hull of points on the solvability set boundary. While it is unlikely that global convexity can be established for realistic power cases, even semi-local convexity covering the realistic operating regions may be extremely helpful in development of next generation of power system decision-support tools.

\section*{Acknowledgment}
The authors would like to thank J.~Bialek for reading the manuscript.

\ifCLASSOPTIONcaptionsoff
  \newpage
\fi



%

\end{document}